# An efficient and easy-to-extend Matlab code of the Moving Morphable Component (MMC) method for three-dimensional topology optimization

Zongliang Du · Tianchen Cui · Chang Liu · Weisheng Zhang · Yilin Guo · Xu Guo*

**Abstract** Explicit topology optimization methods have received ever-increasing interest in recent years. In particular, a 188-line Matlab code of the two-dimensional (2D) Moving Morphable Component (MMC)-based topology optimization method was released by Zhang et al. (Struct Multidiscip Optim 53(6):1243-1260, 2016). The present work aims to propose an efficient and easy-to-extend 256-line Matlab code of the MMC method for three-dimensional (3D) topology optimization implementing some new numerical techniques. To be specific, by virtue of the function aggregation technique, accurate sensitivity analysis, which is also easy-to-extend to other problems, is achieved. Besides, based on an efficient identification algorithm for load transmission path, the degrees of freedoms (DOFs) not belonging to the load transmission path are removed in finite element analysis (FEA), which significantly accelerates the optimization process. As a result, compared to the corresponding 188-line 2D code, the performance of the optimization results, the computational efficiency of FEA, and the convergence rate and the robustness of optimization process are greatly improved. For the sake of completeness, a refined 218-line Matlab code implementing the 2D-MMC method is also provided.

State Key Laboratory of Structural Analysis for Industrial Equipment, Department of Engineering Mechanics, Dalian University of Technology, Dalian 116023, China
Ningbo Institute of Dalian University of Technology, Ningbo 315016, China
E-mail: guoxu@dlut.edu.cn (Xu Guo)

# 1 Introduction

Based on the solid isotropic material with penalization (SIMP) approach, Sigmund developed the celebrated top99 Matlab code (Sigmund, 2001), which greatly promoted the theoretical development and engineering applications of structural topology optimization both in academia and industry. While retaining its compactness, significant speedup was achieved by its upgraded versions top88 and the most recent top99neo (Andreassen et al., 2011; Ferrari and Sigmund 2020). Inspired by these exceptional works, recent years witnessed the emergency of a series of educational articles releasing self-contained compact Matlab codes for disseminating various topology optimization methods, e.g., level set method (Challis, 2010; Wei et al., 2018), evolutionary structural optimization method (Huang and Xie, 2010), Moving Morphable Component (MMC) method (Zhang et al., 2016b), geometry projection method (Smith and Norato, 2020).

Although Matlab provides abundant functions for scientific computing, their efficiency would decrease dramatically when large scale topology optimization problems are solved. Therefore, most of the released Matlab codes are focused on 2D topology optimization problems. Due to the significant computational cost of finite element analysis (FEA) in 3D problems, most existing released Matlab codes are only suitable for medium-sized (with $10^5$ to $10^6$ DOFs) or smaller 3D structural topology optimization problems on a laptop (Ferrari and Sigmund,2020; Liu and Tovar, 2014; Smith and Norato, 2020). To alleviate this issue, the key point is to improve the efficiency of FEA in 3D topology optimization, which constitutes one of the motivations of present work.

In the MMC method, a structure is assumed to be composed by a number of solid components whose sizes, shapes and locations, orientations are described by a set of parameters explicitly (Guo et al., 2014). As illustrated by Fig. 1, by optimizing those geometric parameters, the components could be driven to move, deform, overlap with each other and even disappear to form an optimized structure. Therefore, besides enjoying a crisp structural boundary, in the MMC method, the layout/topology of an optimized structure is described geometrically and explicitly. As a result, relatively little postprocessing work is necessary before importing the obtained design into CAD/CAE system for subsequent analysis/optimization procedures. In addition, in the MMC-based solution framework, feature control and manufacturing constraints can be handled relatively easily by only setting some upper/lower bounds of design variables (Zhang et al., 2016a; Guo et al., 2017; Niu and Wadbro, 2019). Another advantage of the MMC method is its potential to reduce the number of design variables. This makes the MMC method suitable for solving complex topology optimization problems where analytical sensitivities are difficult or even impossible to obtain (e.g., post-buckling problems (Xue et al., 2017), topological material design (Du et al., 2020; Luo et al., 2021), and crashworthiness optimization problems (Raponi et al., 2019), ect). Based on the above reasons, the 188-line Matlab code of the 2D-MMC method (Zhang et al., 2016b) has received considerable interest since its publication.

Although the 3D-MMC method has been developed (Zhang et al., 2017b), a direct extension of the MMC188 Matlab code to 3D case would suffer from the following problems: (1) the analytical sensitivities are not available and this would sacrifice the convergence rate as well as the optimality of the optimized designs; (2) in spite of reducing the number of design variables, how to save the computational cost of FEA, which constitutes the bottleneck for solving 3D topology optimization problems, has not been addressed yet. Although those issues have been partially settled recently (Zhang et al.,

2017; Liu et al., 2018; Xue et al., 2019), the corresponding numerical implementations have not been published in a compact open code.

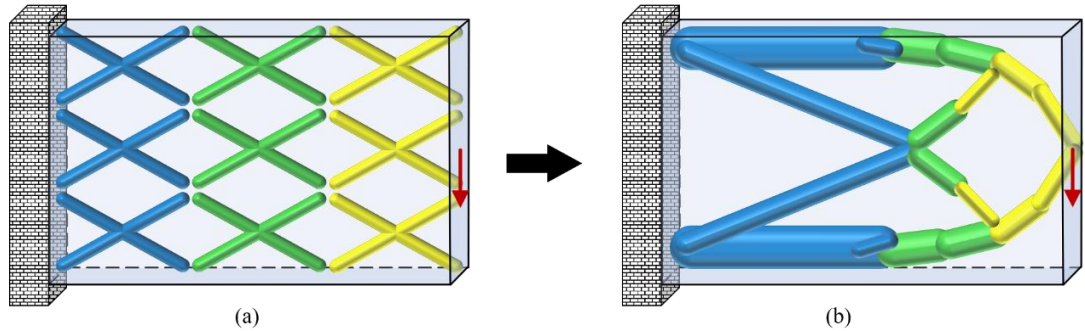

**Fig. 1** An illustration of the 3D-MMC method: (a) initial design, (b) optimized design.

This educational article aims to present a compact and efficient Matlab code, which implements some new effective features of the 3D-MMC method. Through utilizing the Kreisselmeier and Steinhauser function (K-S function, Kreisselmeier and Steinhauser (1980)), the maximum operation for generating the global topology description function is regularized and therefore analytical sensitivities can be obtained to improve the optimality of the optimized designs. Furthermore, a component-based identification algorithm for load transmission path is developed in order to only retain the DOFs pertaining to solid elements in FEA (Zhang et al., 2017a; Xue et al., 2019). It is shown that this treatment is very effective for improving the solution efficiency of 3D topology optimization problems. In addition, some other improvements such as the modifications of the topology description function of components and the introduction of active components and design variables sets are also developed to achieve a robust iteration process with rapid convergence rate.

The rest of the paper is organized as follows. The theoretical aspect of the 3D-MMC method is briefly outlined in Section 2 for the sake of completeness. Then the 256-line Matlab code implementing the 3D-MMC method are explained in detail in Section 3. The numerical performance of the released codes is illustrated in Section 4, by solving some benchmark minimum compliance and compliant mechanism design examples. Finally, the advantages of the provided codes and some possible directions for further improvement are discussed in Section 6.

**2 Moving Morphable Component method for 3D topology optimization problems**

The 3D-MMC method was proposed by (Zhang et al., 2017b). Similar as its two-dimensional counterpart (Guo et al., 2014, 2016; Zhang et al., 2016b), in this approach, a structure is also assumed to be constituted by a number of components whose geometries are expressed by a set of parameters explicitly. As illustrated by Fig. 1, the optimized structure can be obtained by updating the locations, orientations, shapes and sizes of those components through an iteration process. In the present work, to further improve the numerical performance, some modifications to the original 3D-MMC method are made and will be introduced in detail in the following text.

2.1 Topology description of 3D structures using MMCs

In the MMC approach (Guo et al., 2014, 2016; Zhang et al., 2016b, 2017b), the solid structural region occupied by the $i$-th component (i.e., $\Omega^i$) is identified by its topology description function (TDF) $\phi^i$ in the following way:

$$\begin{cases} \phi^i(x) > 0, & \text{if } x \in \Omega^i \cap D \\ \phi^i(x) = 0, & \text{if } x \in \partial\Omega^i \cap D \\ \phi^i(x) < 0, & \text{if } x \in D \setminus (\Omega^i \cup \partial\Omega^i) \end{cases} \quad (1)$$

where D is the design domain.

For the sake of easy implementation, in the present work, all 3D components are restricted to be cuboids, although components with more complicated shapes can also be constructed and used for optimization (Zhang et al., 2017b). Numerical experience showed that this treatment will not sacrifice the performance of optimized solution significantly, as long as the number of components is not too small. By locating a local Cartesian coordinate system $o - x'y'z'$ at the center of the $i$-th component, the corresponding TDF can be constructed approximately as

$$\phi^i = 1 - \left[\left(\frac{x'}{L_1^i}\right)^p + \left(\frac{y'}{L_2^i}\right)^p + \left(\frac{z'}{L_3^i}\right)^p\right]^{1/p} \quad (2)$$

where $L_1^i, L_2^i, L_3^i$ denote the half length, width and height of the component, respectively, $p$ is a relatively large even number (e.g., $p = 6$). Notably, the power $1/p$ in Eq. (2) is introduced to guarantee that the values of the TDF and its partial derivatives do not vary too rapidly.

Denoting the rotation angles between the local and global coordinate systems as $\gamma, \beta, \alpha$, respectively, we have the following relationship for coordinate transformation (see Fig. 2 for reference):

$$\begin{pmatrix} x' \\ y' \\ z' \end{pmatrix} = \begin{bmatrix} \cos\beta^i\cos\gamma^i & \cos\beta^i\sin\gamma^i & -\sin\beta^i \\ \sin\alpha^i\sin\beta^i\cos\gamma^i - \cos\alpha^i\sin\gamma^i & \sin\alpha^i\sin\beta^i\sin\gamma^i + \cos\alpha^i\cos\gamma^i & \sin\alpha^i\cos\beta^i \\ \cos\alpha^i\sin\beta^i\cos\gamma^i + \sin\alpha^i\sin\gamma^i & \cos\alpha^i\sin\beta^i\sin\gamma^i - \sin\alpha^i\cos\gamma^i & \cos\alpha^i\cos\beta^i \end{bmatrix} \begin{pmatrix} x - x_0^i \\ y - y_0^i \\ z - z_0^i \end{pmatrix} \quad (3)$$

where $x_0^i, y_0^i$ and $z_0^i$ denote the coordinates of the center of the $i$-th component in the global Cartesian coordinate system $O - xyz$.

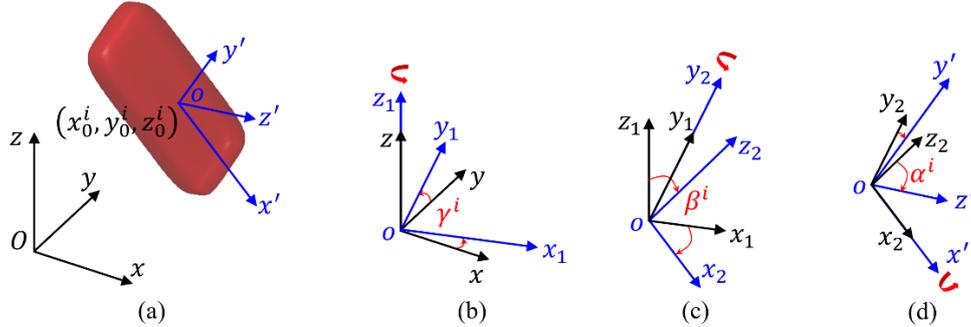

**Fig. 2** Description of a component in local and global Cartesian coordinate systems: (a) an overview, (b-d) rotating the global coordinate system to the local coordinate system.

Combining Eq. (2) and Eq. (3), it yields that the TDF of the $i$-th component in global coordinate system can be described by only 9 geometric parameters explicitly, i.e., $\mathbf{d}^i = \left(x_0^i, y_0^i, z_0^i, L_1^i, L_2^i, L_3^i, \alpha^i, \beta^i, \gamma^i\right)^\top$. Once $\phi^i$ is determined, the TDF of the entire structure can be expressed as $\phi^s = \max(\phi^1, \ldots, \phi^n)$ with $n$ denoting the total number of components introduced in the design domain (Guo et al., 2014; Zhang et al., 2016b, 2017b). To preserve the differentiability of $\phi^s$ necessary for obtaining analytical sensitivities, K-S function[1] is adopted here to approximate the max operation (Liu et al., 2018), i.e.,

$$\phi^s \approx \left(\ln\left(\sum_{i=1}^{n} \exp(\lambda\phi^i)\right)\right)/\lambda \tag{4}$$

where $\lambda$ is a large positive number, e.g., $\lambda = 100$.

As indicated in Eq. (1), the structural domain is actually the region $\Omega = \{x | x \in D, \phi^s(x) > 0\}$ occupied by the solid material. Therefore, the structural layout (material distribution) in the MMC approach is uniquely determined by the vector of the design variables $\boldsymbol{d} = ((\boldsymbol{d}^1)^\top, \dots, (\boldsymbol{d}^n)^\top)^\top$.

*Remark 1* For optimized structures described by tens to hundreds of components, the number of design variables in MMC approach would be only hundreds to thousands even in three-dimensional case. In addition, the design variables are all explicit geometric parameters, this would be very helpful to convert the optimized results to CAD/CAE models for further treatment.

2.2 Problem formulation

For minimum compliance design problems, the mathematical formulation under finite element discretization is expressed as:

$$\begin{aligned}
&\text{find } \boldsymbol{d} = ((\boldsymbol{d}^1)^\top, \dots, (\boldsymbol{d}^n)^\top)^\top \in \mathcal{U}^d \\
&\min f = \boldsymbol{F}^\top \boldsymbol{U} \\
&\text{s.t. } \boldsymbol{KU} = \boldsymbol{F} \\
&\quad g = V/|D| - v \leq 0
\end{aligned} \tag{5}$$

where $\boldsymbol{K}, \boldsymbol{F}$ and $\boldsymbol{U}$ are the global stiffness matrix, nodal force vector and nodal displacement vector; $V, |D|, v$ are the volume of the optimized structure, volume of the design domain, upper bound of allowable volume fraction and $\mathcal{U}^d$ denotes the admissible set of design variables, respectively.

2.3 Finite element analysis (FEA) without redundant DOFs

In the present work, fixed Eulerian mesh and the ersatz material model are adopted for FEA. Under this treatment, the Young's modulus $E_e$ and density $\rho$ of the $e$-th finite element are calculated as

$$E_e = \rho_e E \text{ and } \rho_e = \frac{1}{8}\sum_{i=1}^{8} H_\epsilon^\alpha(\phi_{e,i}^s) \tag{6}$$

where $E$ is the Young's modulus of the solid material and $\phi_{e,i}^s$ is the value of the global TDF $\phi^s$ at the $i$-th node of the $e$-th element. The symbol $H_\epsilon^\alpha$ denotes the smoothed Heaviside function (Zhang et al., 2016b).

$$H_\epsilon^\alpha(x) = \begin{cases} 1, & \text{if } x > \epsilon \\ \frac{3(1-\alpha)}{4}\left(\frac{x}{\epsilon} - \frac{x^3}{3\epsilon^3}\right) + \frac{1+\alpha}{2}, & \text{if } |x| \leq \epsilon \\ \alpha, & \text{otherwise} \end{cases} \tag{7}$$

where $\epsilon$ a parameter controlling the size of the transition zone of $H_\epsilon^\alpha$ and $\alpha$ is a small positive parameter (e.g., $\alpha = 10^{-3}$) introduced to avoid the possible singularity of the global stiffness matrix associated with a specified design. It should be point out that, the smoothed Heaviside function here is a S-shape mapping function, which can effectively pronounce the sensitivity of cut elements with intermediate densities and push the evolution of structural boundaries. Interested readers is referred to Yoon and Kim (2003) for more details of the role of S-shape mapping functions in structural topology optimization.

---

[1] It should be noted that, alternatively, the $p$-norm function is another choice for approximating the max operation and analytical sensitivities could also be obtained accordingly (Picelli et al., 2018).

Once the Young's modulus of each element is determined, FEA can be performed straightforwardly. Furthermore, as long as a load transmission path is formed, the redundant DOFs only associated with the elements where $\rho = \alpha$ can be totally discarded in FEA without sacrificing the accuracy (Zhang et al., 2017a, 2018; Xue et al., 2019). This treatment can not only save the computational time but also eliminate the numerical instabilities caused by weak elements. The implementation of this redundant DOFs removal technique and its effectiveness will be discussed and illustrated in the illustration example section.

2.4 Sensitivity analysis

For the concerned objective and constraint functions $f$ (structural compliance) and $g$ (volume of the structure) in Eq. (5), the corresponding sensitivities with respect to a specific design variable $d_j^i$ (the $j$-th parameter associated with the $i$-th component) can be calculated as:

$$\frac{\partial f}{\partial d_j^i} = \sum_{e=1}^{N} \frac{\partial f}{\partial \rho_e}\frac{\partial \rho_e}{\partial d_j^i} = -\sum_{e=1}^{N} \boldsymbol{u}_e^\top \mathbf{k}_e^0 \boldsymbol{u}_e \frac{\partial \rho_e}{\partial \phi^s}\frac{\partial \phi^s}{\partial \phi^i}\frac{\partial \phi^i}{\partial d_j^i} \tag{8}$$

and

$$\frac{\partial g}{\partial d_j^i} = \sum_{e=1}^{N} \frac{\partial g}{\partial \rho_e}\frac{\partial \rho_e}{\partial d_j^i} = \sum_{e=1}^{N} \frac{V_e}{|\mathrm{D}|}\frac{\partial \rho_e}{\partial \phi^s}\frac{\partial \phi^s}{\partial \phi^i}\frac{\partial \phi^i}{\partial d_j^i} \tag{9}$$

where $u_e, V_e$ are the corresponding nodal displacement vector and volume of the $e$-th element; $\mathbf{k}_e^0$ is the stiffness matrix of an element composed of the solid material. According to Eq. (4) and Eq. (6), it yields that

$$\frac{\partial \rho_e}{\partial \phi^s}\frac{\partial \phi^s}{\partial \phi^i} = \frac{1}{8}\sum_{j=1}^{8} H_\epsilon^{\alpha'}(\phi_{e,j}^s) \frac{\exp(\lambda \phi_{e,j}^i)}{\sum_{k=1}^{n}\exp(\lambda \phi_{e,j}^k)} \tag{10}$$

The detailed expressions of $\frac{\partial \phi^i}{\partial d_j^i}$ for the adopted cuboid MMC are presented in Appendix 1.

*Remark 2* The above sensitivity results can be easily generalized to the cases where general objective/constraint functions are involved. This is because for a specified kind of component, the term $\partial \rho_e/\partial d_j^i$ is the same and independent on the concerned objective/constraint functions. In addition, the other term $\partial f/\partial \rho_e$ is nothing but the pre-filtered sensitivity in SIMP method when the penalization power is set to be 1. Therefore, all sensitivity results obtained in SIMP method can be used directly in the present MMC-based approach.

**3 The implementations of the 256-line 3D-MMC Matlab code**

In this section, taking the cantilever beam illustrated in Fig. 3(a) for example, the MATLAB code (see Appendix 4) implementing the 3D MMC-based structural topology optimization approach is explained in detail.

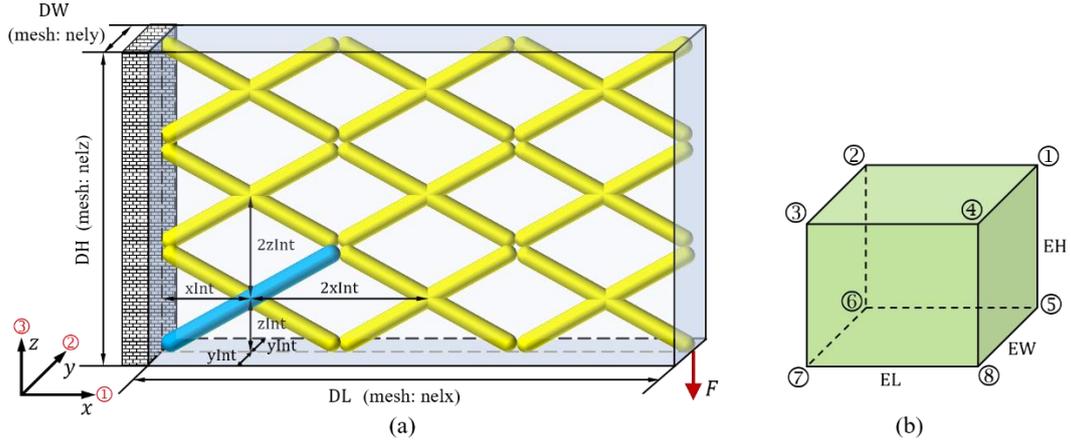

**Fig. 3** (a) An illustration of the initial setting of the 3D cantilever beam problem, (b) geometry description of a eight-node brick element.

3.1 Main function

In this part, some parameters are initialized. As shown in Fig. 3, the cuboid design domain is assumed to have a size of DL × DW × DH and discretized by nelx × nely × nelz uniform eight-node brick elements. The maximum volume fraction is volfrac. The initial design can be generated from a representative component (highlighted in blue in Fig. 3(a)). Its center coordinate is denoted as (xInt, yInt, zInt), its half-length, half-width and half-height are denoted as L1, L2, L3, its three rotation angles are denoted as alp, bet, gam, respectively.

One can activate MMC3d.m by executing the following two command lines:

vInt = [L1 L2 L3 alp bet gam];

MMC3d(DL, DW, DH, nelx, nely, nelz, xInt, yInt, zInt, vInt, volfrac);

In the following, the MMC3d.m is explained in detail.

3.2 Parameter setting (lines 2-7)

In these lines, besides defining the isotropic elastic constants of the solid material, the values of the parameters dgt0 (the number of reserved significant digits of sensitivities to alleviate the accumulation of truncation errors) and scl (the factor used to keep a balance between the scaled values of the objective and the constraint functions) are also set. The parameters p and lmd correspond to the quantities $p$ and $\lambda$ in Eq. (2) and Eq. (4), respectively. Besides, the relative variation of the values of the objective function in the latest consecutive five iterations (i.e., objVr5) is setting as:

$$objVr5(k) = \begin{cases} 1.0, \text{ if } k < 5 \\ objVr5(k-1), \text{ if } k \geq 5 \text{ and } V_{er} > 10^{-4} \\ \left|\frac{max(max(l) - \sum_{k-4}^{k} obj(k)/5\ |}{\sum_{k-4}^{k} obj(k)/5}\right|, \text{ else} \end{cases} \quad (11)$$

where $V_{er} = \left(\frac{V(\boldsymbol{D})}{|\boldsymbol{D}|} - v\right)/v$ with $V(\boldsymbol{D})$ denoting the volume of the current design.

3.3 Setting of FE discretization (lines 8-32)

The design domain is discretized by nEle finite elements, nNod nodes and nDof degree of freedoms. The smallest size of the uniform finite elements is minSz. Both the FE nodes and elements are numbered from the left-bottom corner, first along the $x$ direction, then along $y$ direction, and finally counting layer-by-layer along $z$ direction as illustrated by Fig. 3(a). Accordingly, nNodfc is the number of nodes in each $x-y$ layer. The symbols alpha and epsilon denote the values of $\alpha$ and $\epsilon$ in the smoothed Heaviside function of Eq. (7).

By calling the Ke_tril function in lines 144-177, lines 15-17 produce the elemental stiffness matrix of the eight-node brick element as illustrated in Fig. 3(b). Following the treatment in Ferrari and Sigmund (2020a), the matrices in lines 18-22 are stored in int32 type. Besides, the operations in lines 23-29 prepare for assembling the global stiffness matrix by taking the symmetry property into consideration[2].

In line 30, the structured array LSgrid stores the coordinates of each nodes used for the calculation of the TDFs. The symbol volNod in line 32 is the weight vector of each node for calculating the structural volume.

3.4 Loads, displacement boundary conditions (lines 33-40)

The indices of fixed nodes (fixNd), fixed DOFs (fixDof), elements containing fixed nodes (fixEle), free DOFs (freeDof), loading nodes (loadNd), loading DOFs (loadDof), elements containing loading nodes (loadEle) and the vector of the external load (F) are generated.

3.5 Initial setting of components (lines 41-60)

In lines 42-46, the vectors of the initial central coordinates of the components are saved as x0, y0, z0. Similarly, in lines 48-53, the initial half-lengths, half-widths, half-heights and the inclined angles of the components are stored in $1 \times N$ row vectors.

In lines 57-58, all the components and the corresponding design variables are set as active in the first step and saved in actComp and actDsvb, respectively.

In the present implementation, the non-designable domains are treated as fixed components, and nNd and PhiNd are used to describe the numbers of the non-designable domains and vectors of their nodal TDF values (line 59).

In line 60, the TDF values of all components are initialized in allPhi.

3.6 Setting of MMA (lines 61-68)

For minimum compliance design problems subject to a volume constraint, the MMA parameters are initialized in this part.

3.7 Optimization loop (lines 69-141)

*3.7.1 Generating TDFs and their derivatives (lines 71-82)*

In each iteration, the nodal values of the derivatives of TDFs for each component with respect to their own design variables are initialized and stored in allPhiDrv. Then for the $i$-th component belonging to the current active components set actComp, the nodal values of TDF and the corresponding sensitivities are calculated by the calc_Phi function (lines 178-213) and stored as:

$$\text{allPhi} = [\phi^1 \dots \phi^i \dots \phi^{N+nNd}]_{nNod \times (N+nNd)} \tag{12}$$

$$\text{allPhiDrv} = \left[\frac{\partial \phi^1}{\partial d_1^1} \dots \frac{\partial \phi^i}{\partial d_1^i} \dots \frac{\partial \phi^i}{\partial d_9^i} \dots \frac{\partial \phi^N}{\partial d_9^N}\right]_{nNod \times 9N} \tag{13}$$

---

[2] The fsparse routine developed by Engblom and Lukarski (2016) is required (available at https://github.com/stefanengblom/stenglib).

To improve the convergence rate, in lines 206-212, tiny component (either with $\min(L1; L2; L3)/minSz \leq 0.1$ or all its absolute TDF values no smaller than $\epsilon$), would be set as null and removed from the active set. Otherwise, $\partial \phi^i / \partial d_j^i$ are calculated according to lines 193-205.

Subsequently, from lines 77 to 79, the TDFs' nodal values of the active components and the non-design solid domains are collected in allPhiAct and further used to calculate the nodal value vector of the global TDF (Phimax) implementing Eq. (4)[3]. Then the sensitivities of the global TDF vector with respect to active design variables are obtained as PhimaxDrvAct in line 82.

*3.7.2 Plotting current design (lines 83-88)*

The current design is plotted in this subsection, and one could also only plot the final design for saving time.

*3.7.3 Finite element analysis (lines 89-118)*

For the volume constraint, the nodal density vector (H) of current design is calculated using the smoothed Heaviside function (Heaviside) in lines 247-252. Then the elemental density vector is obtained at line 91 as Eq. (6). For FEA, as described in Subsection 2.3, in order to apply the redundant DOFs removal technique, a load transmission path should be identified.

To this end, an efficient component-based identification algorithm for load transmission path, as explained by Fig. 4, is implemented in the function srch_ldpth in lines 214-246. For the case where a load transmission path exists, strct returns 1 and loadPth at line 94 returns the index numbers of related components[4] ; otherwise, strct returns 0.

In the MMC method, a load transmission path is composed by components continuously connecting the external load and fixed boundary. As a basic requirement of the existence of a load transmission path, therefore, in the subroutine srch_ldpth, it first makes sure all the loading elements and at least one of the fixed elements are non-void (line 223). Based on that, utilizing the connectivity relation between active components, one can find the set of components (the blue components in Fig. 4(d)) which can be connected to the loading components by step, as illustrated by Fig. 4(a)-(d). If some component in this set is also constrained by the fixed boundary, the path exists and this set of components constitute the load transmission path.

For the case a load transmission path existing, as indicated in lines 96-101, the elemental density vector for FEA is updated as denSld only based on the components belonging to the load transmission path. Then the DOFs only pertaining to the elements whose density values are smaller than alpha + eps are treated as redundant. Therefore, only the retained free DOFs, i.e., freedofLft, are involved in the FEA at line 110. Notably, the global stiffness matrix is efficiently assembled in line 108, following the treatment in Ferrari and Sigmund (2020). Notably, to avoid the possible singularity caused by a disconnected density field resulting from some peculiar configurations of thin components (see Fig. 5 for reference), besides the introduction of inactive components, a small value (denoted by eps) is added to the diagonal elements of the global stiffness matrix (line 109).

For an intermediate design without a load transmission path, in lines 113-115, the FEA is implemented as usual (Andreassen et al., 2011; Zhang et al., 2016b; Ferrari and Sigmund, 2020).

---

[3] Phimax is regularized with a minimum value of $10^{-3}$ to avoid possible singularity.
[4] The non-design solid domains also belong to the load transmission path.

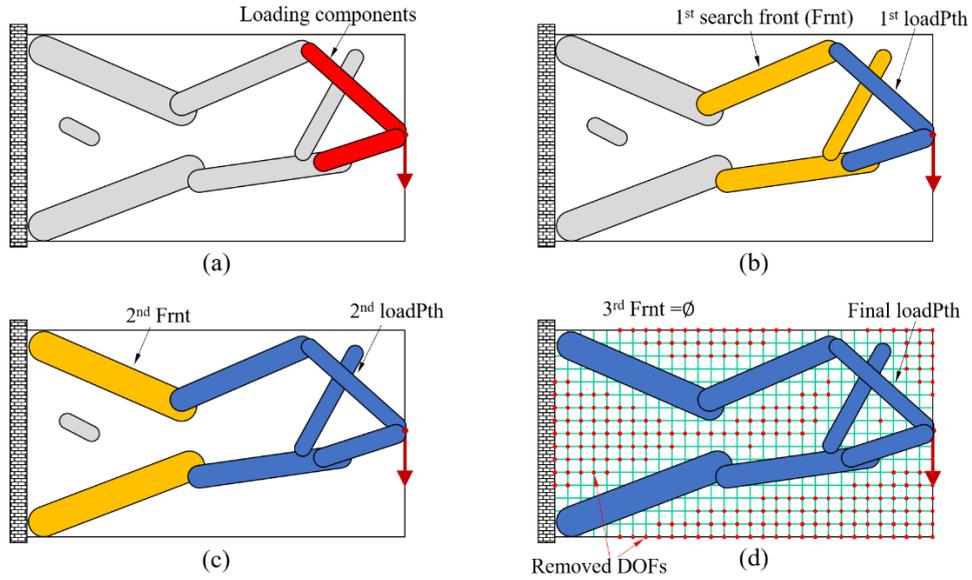

**Fig. 4** Illustration of the identification algorithm for load transmission path. (a) Find the loading components (highlighted in red), (b) collect their index numbers as the elements of the first load path set (loadPth in line 232, and related components are highlighted in dark blue), and according to the connectivity among active components (cnnt in line 228), collect the index numbers of components (highlighted in yellow) connecting to components in the first search front set (Frnt in line 233); (c) combine the current Frnt into the loadPth (line 238), and then collect the index numbers of the components connecting to components identified by the current Frnt as the Frnt of the next step in line 242, (d) when the current Frnt is empty, the identification process for load transmission path terminates (line 244). If a load transmission path is identified, the corresponding redundant DOFs (related to nodes colored in red) can be eliminated from the FEA.

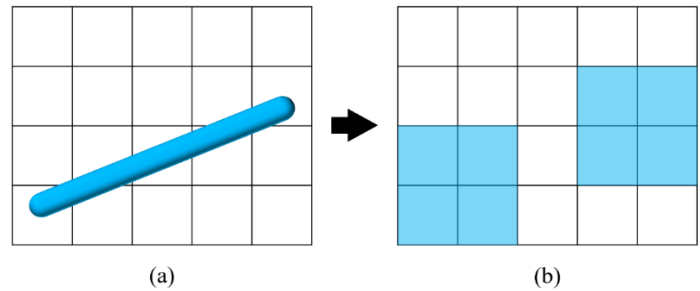

**Fig. 5** (a) A thin component in the design domain, (b) the corresponding disconnected density field obtained using the adopted ersatz material model.

*3.7.4 Sensitivity analysis (lines 119-130)*

In line 122, the nodal values of $\partial H_\epsilon^\alpha/\partial \phi^s$ are stored as delta_H. Since delta_H and PhimaxDrvAct are both stored in the nodal manner, the quantity $\boldsymbol{u}_e^\top \boldsymbol{k}_e^0 \boldsymbol{u}_e$ in Eq. (8) is also stored in its nodal values (engyNod) in line 125. Under this circumstance, the formula of sensitivities of structural compliance (df0dx) and volume constraint (dfdx) in Eq. (8) and Eq. (9) are executed in their nodal values in lines 126-127 as:

$$\frac{\partial f}{\partial d_j^i} = -\sum_{m=1}^{nNod} E_m \frac{\partial H_{\epsilon,m}^{\alpha}}{\partial \phi_m^s} \frac{\partial \phi_m^s}{\partial \phi_m^i} \frac{\partial \phi_m^i}{\partial d_j^i} \tag{14}$$

and

$$\frac{\partial g}{\partial d_j^i} = \sum_{m=1}^{nNod} \frac{V_e}{|D|} W_m \frac{\partial H_{\epsilon,m}^{\alpha}}{\partial \phi_m^s} \frac{\partial \phi_m^s}{\partial \phi_m^i} \frac{\partial \phi_m^i}{\partial d_j^i} \tag{15}$$

where the quantities $E_m$, $W_m$, $\frac{\partial H_{\epsilon,m}^{\alpha}}{\partial \phi_m^s}$ (vector) and $\frac{\partial \phi_m^s}{\partial \phi_m^i} \frac{\partial \phi_m^i}{\partial d_j^i}$ (matrix) are denoted as engyNod, volNod, delta_H and PhimaxDrvAct in the program, respectively.

Finally, lines 128-130 round the sensitivities by only keeping dgt0 significant digits to alleviate truncation errors.

*3.7.5 Updating design variables (lines 131-140)*

The method of moving asymptotes (MMA, Svanberg(1987)) is used as the optimizer.[5] In the mmasub.m, the corresponding parameters are set as epsimin = 1e − 9; raa0 = 0.01; albefa = 0.8; asyinit = 0.05; asyincr = 1.0; asydecr = 0.8. It is worth noting that p0 = max(df0dx, 0) and q0 = max(−df0dx, 0) in the original mmasub.m should be revised as p0 = max(df0dx′, 0) and q0 = max(−df0dx′, 0) in accordance with the data format adopted in the present implementation.

After updating the design variables, in lines 135-137, the relative variation of the objective function values in latest five consecutive iterations (objVr5) is updated.

**4 Illustration examples**

In this section, four illustrative examples are investigated by modifying the proposed 256-line MMC3d code to demonstrate its features including the influence of different parameters $\epsilon$ on the optimized designs, the accuracy of analytical sensitivities, the computational efficiency of FEA with the redundant DOFs removal technique, and its performances for non-design domains and different objective functions. Besides, a refined 218-line MMC2d code is also verified by the fifth example. All the examples are run on a laptop equipped with an Intel(R) Core(TM) i9-10885H@2.40GHz 2.40GHz CPU, 64.0GB of RAM, and Matlab 2020b under Windows 10.

4.1 Cantilever beam example

The cantilever beam example studied in Hoang and Nguyen-Xuan (2020) is revisited here, as illustrated in Fig. 3(a). The main code is specified as:

vInt = [12 2.5 2.0 0 atan(1) 0];
MMC3d(64, 8, 32, 128, 16, 64, 8.0, 4.0, 8.0, vInt, 0.3);

As illustrated by Fig. 6(a), the initial design is composed of 16 MMCs and the number of design variables is $16 \times 9 = 144$. In addition, the comparison in Fig. 6(b) reveals the analytical sensitivities match well with the corresponding finite difference results (central difference, $\Delta = 10^{-8}$).

The optimization loop is terminated at iteration 171 when $objVr5 = 9.76 \times 10^{-5} < 10^{-4}$, and the iteration history of the objective and constraint function values are plotted in Fig. 7. Some representative intermediate designs and the corresponding structural compliances are presented in Fig. 8(a-d). At iteration 7, the right-top one in Fig. 8(b) is too thin ($L_3/\min(EL, EW, EH) < 0.1$), and is set as inactive. At around the 100-th iteration, the optimized design is evolved to an

---
[5] The MMA code of 2007 version is used (Svanberg, 2007). Available at http://www.smoptit.se/

*I*-shaped beam which is well-known for its excellent bending moment resistance capacity. After some minor improvement in structural configuration such as the smoothness of the core plate and the gradient of the bottom illustrated in Fig. 8(e), compared to the reference optimized design using extruded geometric components (Fig. 8(f), Hoang and Nguyen-Xuan (2020)), the value of the objective function is decreased by about 10.4% (20.50 vs 18.37).

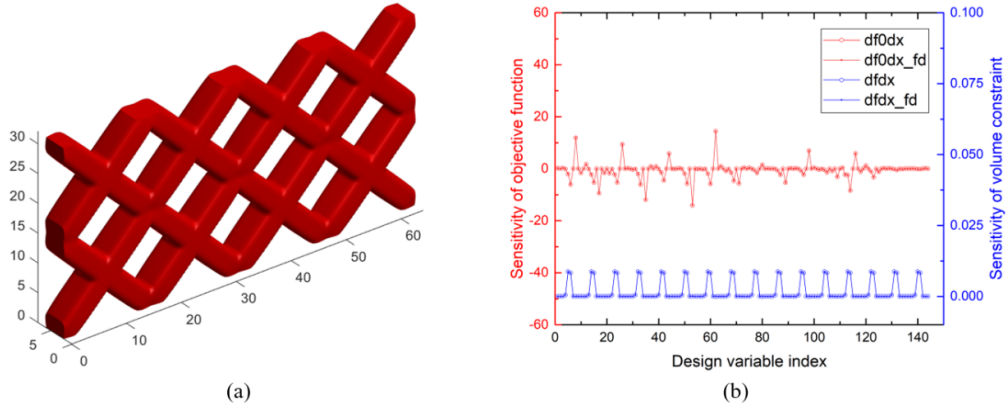

(a)  (b)

**Fig. 6** (a) The initial design of the 3D cantilever beam example, (b) a perfect match between the analytical and finite difference sensitivities.

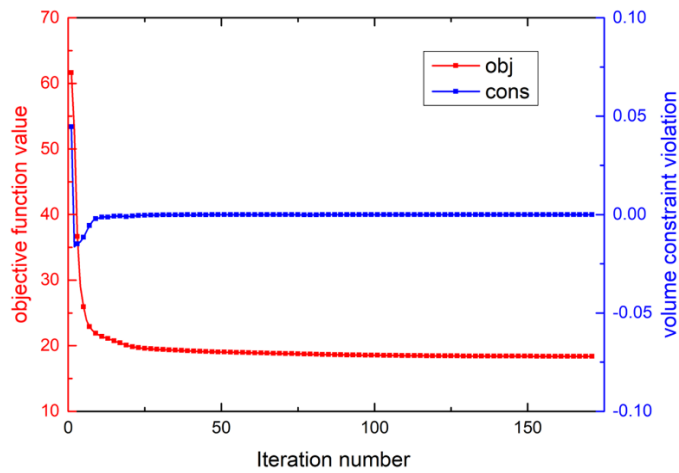

**Fig. 7** The iteration history of the 3D cantilever beam example.

With other problem setting unchanged, this problem is also solved with $\epsilon = 0.1$ and $\epsilon = 0.5$ respectively. For comparison, the iteration histories $\epsilon = 0.1, 0.25, 0.5$ are illustrated in Fig. 9. As indicated in Table 1, the corresponding iteration numbers for achieving the same convergence requirement ($objVr < 10^{-4}$) are 297, 171 and 115, respectively. This can be partially explained by the fact that as $\epsilon$ gets larger, there would be more elements contributing to the sensitivities of the objective and constraint functions and thus the iteration process may converge faster.

*Remark 3* According to Eq. (7), the treatment by setting $|x| \leq \epsilon$ is similar to the narrow band technique often adopted in the level set topology optimization method (Wang et al. 2003). It is worth noting, however, that it is more reasonable to choose $0.1 \leq \epsilon \leq 0.5$ in the MMC method. This is because the MMC-based TDF of a structure has a maximum value of 1 and its nodal values are independent on the absolute size of the structure. For a smaller value of $\epsilon$, the convergence rate would be slow, while for an overly large value of $\epsilon$, there would be too many gray elements in the finite element model.

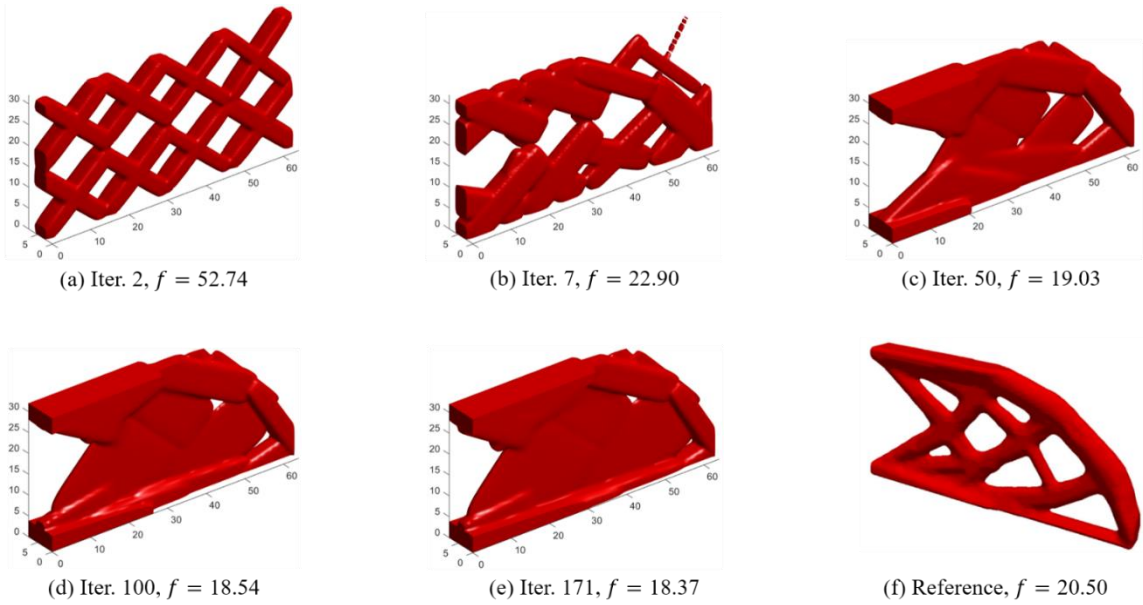

(a) Iter. 2, $f = 52.74$
(b) Iter. 7, $f = 22.90$
(c) Iter. 50, $f = 19.03$
(d) Iter. 100, $f = 18.54$
(e) Iter. 171, $f = 18.37$
(f) Reference, $f = 20.50$

**Fig. 8** (a)-(e) Representative intermediate designs and their compliance values, (f) reference design and its compliance value obtained in Hoang and Nguyen-Xuan (2020).

**Table 1** Comparison of the optimized designs obtained with different values of $\epsilon$ ($\lambda = 100$).

|  | $\epsilon = 0.1$ | $\epsilon = 0.25$ | $\epsilon = 0.5$ |
|---|---|---|---|
| Iteration | 297 | 171 | 115 |
| Compliance | 18.56 | 18.37 | 18.39 |
| Configuration | 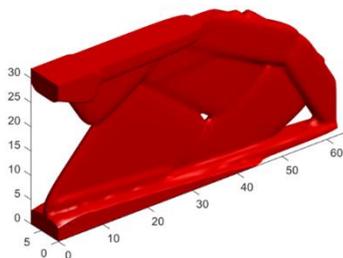 | 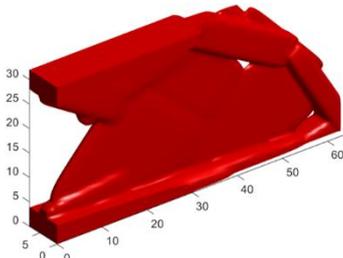 | 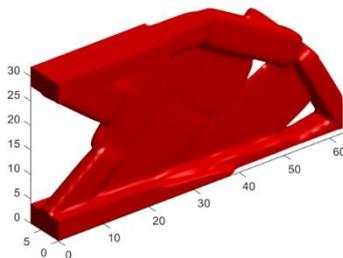 |

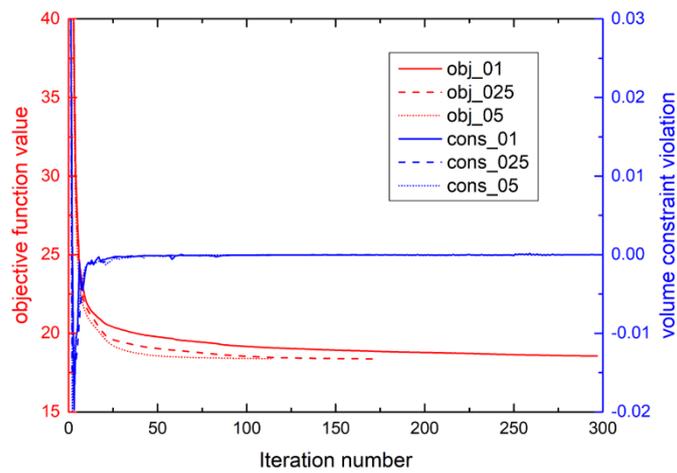

**Fig. 9** The iteration histories of the cantilever beam example with different values of $\epsilon$ ($\lambda = 100$).

As illustrated by Table 2, different values of $\lambda \geq 20$ would slightly affect the optimized structural compliance. Notably, to better approximate the max operation, a larger $\lambda$ is desired, while an overlarge value would lead to the fluctuation of iteration history as show in Fig. 10.

**Table 2** Comparison of the optimized designs obtained with different values of $\lambda$ ($\epsilon = 0.25$).

| | $\epsilon = 0.1$ | $\epsilon = 0.25$ | $\epsilon = 0.5$ |
|---|---|---|---|
| Iteration | 165 | 171 | 272 |
| Compliance | 18.37 | 18.37 | 18.39 |
| Configuration | 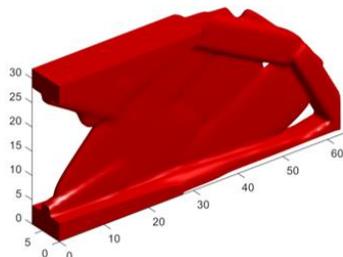 | 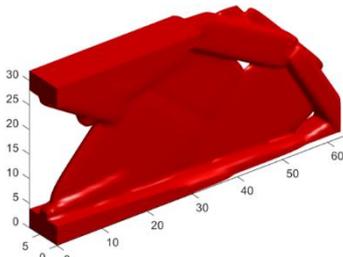 | 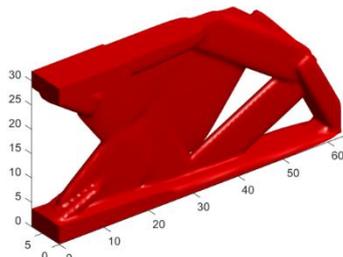 |

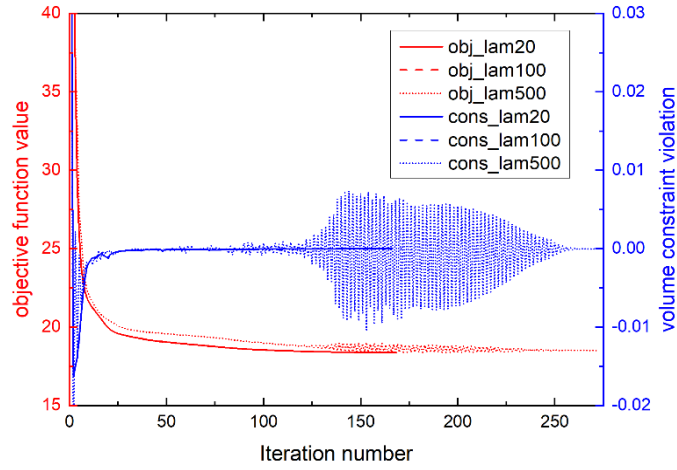

**Fig. 10** The iteration histories of the cantilever beam example with different values of $\lambda$ ($\epsilon = 0.25$).

4.2 MBB beam example

The MBB beam example illustrated in Fig. 11(a) is optimized with a volume constraint of $V/|D| \leq 0.25$. Due to the symmetry of the problem setting, only the light-green part in Fig. 11(a) is optimized. The corresponding main code is:

vInt = [0.7 0.08 0.08 atan(1/2) atan(1/2) atan(-1/2)];

MMC3d(3, 0.5, 1, 60, 10, 20, 0.5, 0.25, 0.25, vInt, 0.25);

In the MMC3d.m, the corresponding modifications are made:

(1) In SEC 2), the $\epsilon$ is modified as 0.2;

(2) All the SEC 3) about loads and displacement boundary conditions from lines 34-40 are updated as follows:

[jN,kN] = meshgrid(1:nely+1,1:nelz+1);  fixNd1 = [(kN-1)*nNodfc+(jN-1)*(nelx+1)+1];

[iN,kN] = meshgrid(1:nelx+1,1:nelz+1);  fixNd2 = [kN*nNodfc-iN+1];

fixNd3 = nelx + 1;

fixDof = [3*fixNd1(:)-2; 3*fixNd2(:)-1; 3*fixNd3(:); 3*fixNd3(:)-1; 3*fixNd3(:)-2];

fixEle = [nelx];

```
freeDof = setdiff([1:nDof],fixDof);
loadNd = nNod - nelx; loadDof = 3*loadNd;
loadEle = [nEle-nelx+1];
F = fsparse(loadDof',1,-1/4,[nDof,1]);
```
(3) In SEC 4), lines 44-46 are modified as
```
x0 = kron(x0,ones(1,4*yn*zn));
y0 = repmat(repmat(y0,1,2*zn),1,2*xn);
z0 = repmat(reshape(repmat(z0,4,1),1,8),1,xn);
```

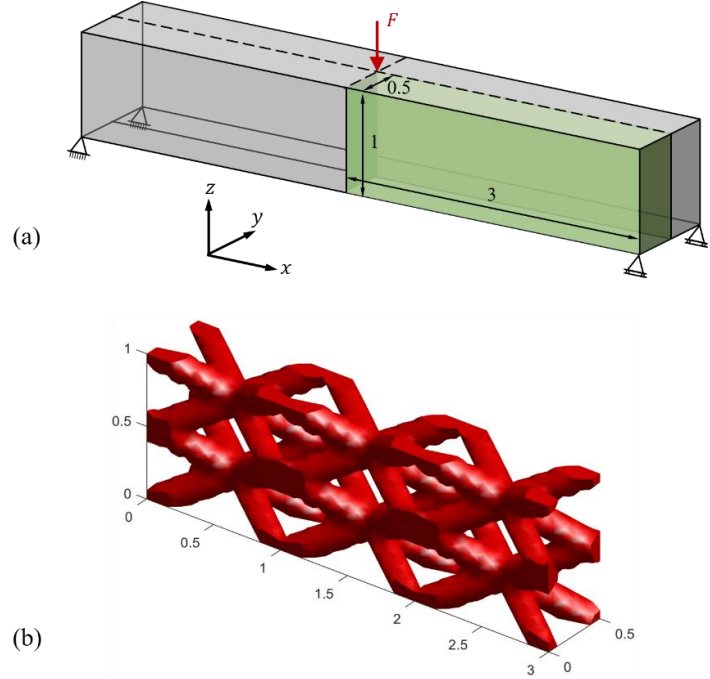

**Fig. 11** (a) An illustration of the MBB example, (b) the initial distribution of the components.

The initial design with 24 components (only $24 \times 9 = 216$ design variables) is illustrated in Fig. 11(b). The optimization loop is terminated at iteration 192 with an optimized structural compliance of 40.76. Notably, all the intermediate structures have a load transmission path, and thus the FEA is implemented with the help of redundant DOFs removal technique at every step of optimization. To test its efficiency, the intermediate structures are also analyzed with full DOFs. The iteration histories of the objective function value, the relative time costs of searching for load transmission path ($t_{srch}$) and FEA with reduced DOFs ($t_{FEr}$) with respect to the cost of FEA with full DOFs ($t_{FE}$) are illustrated in Fig. 12. Obviously, the relative time cost of the load transmission path identification process is very low (always less than 5%). Besides, in most iterations, solving equilibrium equation of reduced DOFs only takes about 35% of the FEA cost when full DOFs are included.

To further study the performance of the present code for different problem sizes, this MBB example is also solved with and without using DOFs removal technique under three types of meshes (i.e., $60 \times 10 \times 20$, $90 \times 15 \times 30$ and $120 \times 20 \times 40$ uniform finite elements), respectively. The corresponding optimized results and average time cost per iteration are listed in Table 3.

It can be observed that the configurations and the structural compliance values of each pair of the optimized designs are actually quite close. As the FE mesh is refined, the holes in the vertical part of the beam are closed and part of the cross section tends to take a I-shape (view along the $x$-axis) since the bending moment is relatively large at the middle part of the beam.

**Table 3** The optimized MBB beams and the time costs of FEA with different mesh sizes.

| Mesh size | Remv. | Iter. | Optimized design | Compliance | Average time per iter. (s) Load path | Average time per iter. (s) FE reduced |
|---|---|---|---|---|---|---|
| 60 × 10 × 20 | with | 192 | 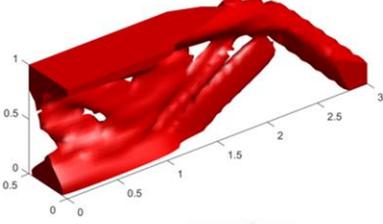 | 40.76 (59.31) | 0.017 (2.2%) | 0.237 (30.2%) |
| 60 × 10 × 20 | without | 500 | 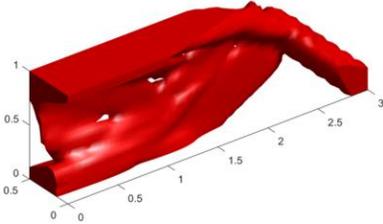 | 40.43 (58.74) | 0.786 | |
| 90 × 15 × 30 | with | 500 | 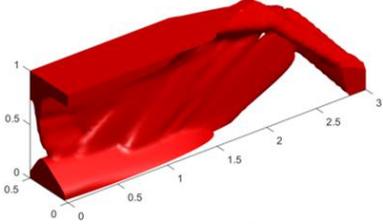 | 49.14 (57.71) | 0.057 (0.6%) | 1.033 (11.7%) |
| 90 × 15 × 30 | without | 500 | 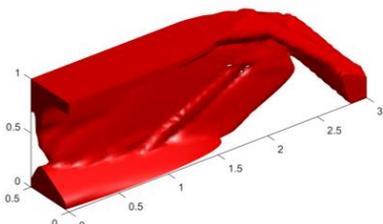 | 49.04 (57.60) | 8.812 | |
| 120 × 20 × 40 | with | 326 | 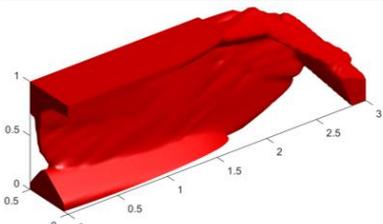 | 57.73 (57.66) | 0.140 (0.2%) | 3.956 (7.0%) |
| 120 × 20 × 40 | without | 394 | 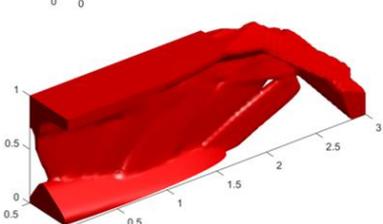 | 57.63 | 56.587 | |

Notably, the compliance values of the optimized designs increase significantly, i.e., from 40.76 to 57.74 by doubling the mesh on each direction. This can be well-understood by the shear locking phenomena associated with the first-order, fully-integrated elements subject to bending deformation (Bathe, 2006). The structural compliances obtained by FEA with full DOFs under a 120 × 20 × 40 mesh are also presented (inside the bracket) for comparison purpose. The provided results indicate when stiffness design of bending-dominated structures is considered, one should take the shear locking effect into consideration seriously.

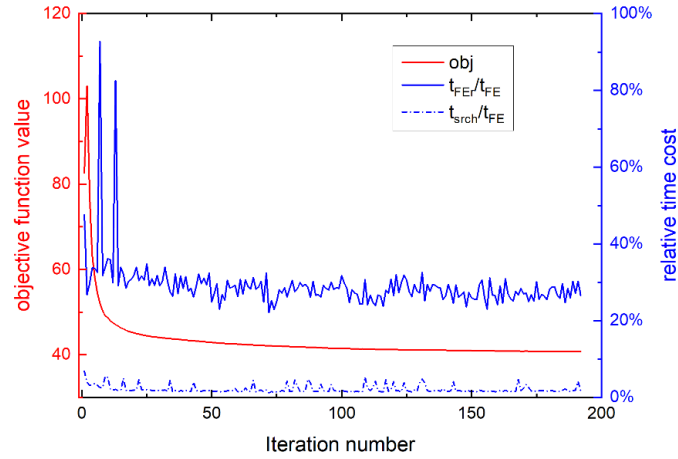

**Fig. 12** The iteration histories of the objective function values, the relative time costs of the FEA with redundant DOFs removal technique (MBB example).

4.3 Torsion beam example with non-design domains

In this example, as illustrated in Fig. 13, a $12 \times 4 \times 4$ beam is fixed at left-hand side and four concentrated loads are applied on the corners of right-hand side. Two faces with width of 0.25 on the left and right are fixed as non-design solid domains. The objective is to minimize the structural compliance with a maximum allowable volume fraction of 15%. Similar to the treatment in Zhang et al. (2017a), the design domain is discretized by $96 \times 32 \times 32$ uniform eight-node brick elements.

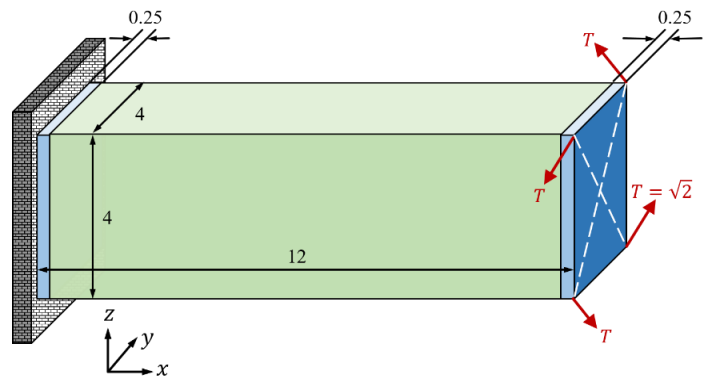

**Fig. 13** An illustration of the torsion beam example.

Accordingly, the main code is revised as:

vInt = [1.2 0.2 0.2 0 asin(0.4) asin(-0.4)];

MMC3d(12, 4, 4, 96, 32, 32, 1.0, 0.5, 0.5, vInt, 0.15);

In the MMC3d.m, the modifications of the problem settings are listed as follows:

(1) In SEC 1, the scale factor (scl) is set as 1000 to make the objective function more consistent with the volume constraint.

(2) In SEC 2, the value of epsilon is modified as 0.5 to improve the convergence rate.

(3) In SEC 3, lines 38-40 describing the load boundary condition are updated as:

loadNd = [nelx+1, nNodfc, nNod-[(nelx+1)*nely,0]];

loadDof = [3*loadNd(1)-[1 0],3*loadNd(2)-[1 0],...

3*loadNd(3)-[1 0],3*loadNd(4)-[1 0]];

loadEle = [nelx,nelx*nely,nEle-[nelx*nely 0]];

F = fsparse(loadDof',1,[1, -1, 1, 1, -1, -1, -1, 1]',[nDof,1]);

(4) In SEC 4, lines 44-46 and lines 52-53 describing the initial distribution of the components are modified as:

x0 = kron(x0,ones(1,yn*zn));

y0 = repmat(repmat(y0,1,zn),1,xn);

z0 = kron(repmat(z0,1,xn),ones(1,yn));

and

tem1 = [1 1 1 1]*vInt(5);    bet = repmat([tem1 -tem1 tem1 -tem1 -tem1 tem1 -tem1 tem1],1,3);

tem2 = [1 -1 1 -1]*vInt(6);   gam = repmat([tem2 tem2 tem2 tem2 -tem2 -tem2 -tem2 -tem2],1,3);

In addition, the command of line 59 about the non-design domains is revised as:

nNd = 2;    PhiNd = -1*ones(nNod,2);

nDnd1 = fixNd(:) + [0:2];   nDnd2 = fixNd(:) + nelx - [0:2];

PhiNd(nDnd1(:),1) = 1;  PhiNd(nDnd2(:),2) = 1;

The initial design is illustrated in Fig. 14(a), where totally 96 components and 864 design variables are introduced. In the 4-th and 5-th iterations (Figs. 14(b,c)), there are 24 tiny components set as inactive, respectively. Thereafter, as shown in Fig. 14(d), the remaining components together form a reticular structure and only 432 design variables are involved. At iteration 15, the reticular structure gets expanded and its configuration and structural compliance value are both close to the optimized design obtained by the original 3D-MMC method (Zhang et al., 2017a,b). Actually, a further optimized design can be obtained at iteration 227 (Fig. 14(f)), with a 30.6% decrease of the structural compliance value (i.e., 1902.86 vs 2742.12). The non-smooth surface can be improved by increasing the mesh resolution as shown later.

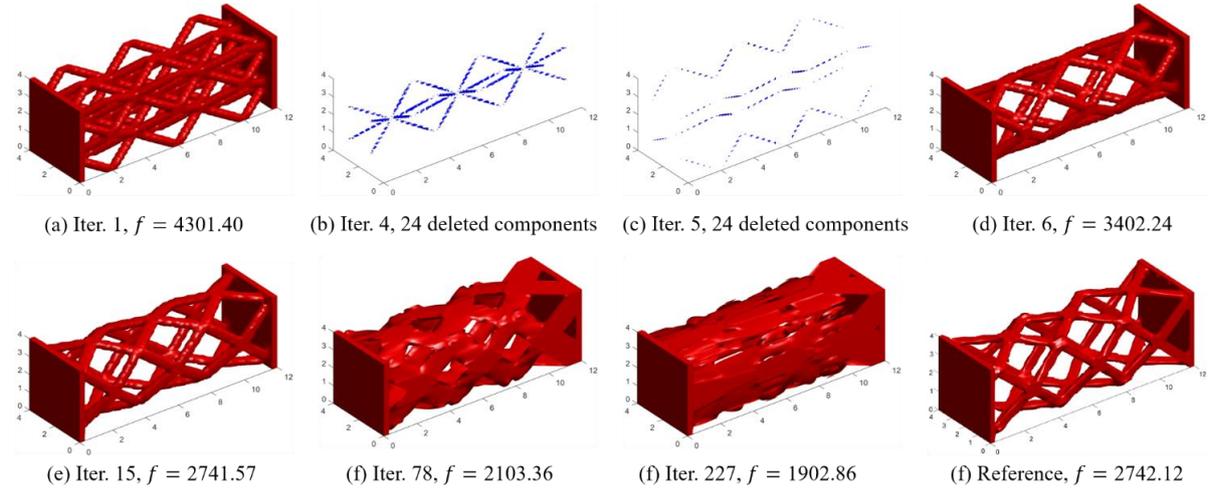

(a) Iter. 1, $f = 4301.40$  (b) Iter. 4, 24 deleted components  (c) Iter. 5, 24 deleted components  (d) Iter. 6, $f = 3402.24$

(e) Iter. 15, $f = 2741.57$  (f) Iter. 78, $f = 2103.36$  (f) Iter. 227, $f = 1902.86$  (f) Reference, $f = 2742.12$

**Fig. 14** (a)-(f) Representative intermediate configurations of the torsion beam example with a $96 \times 32 \times 32$ mesh, (f) reference design and its compliance value obtained in Zhang et al. (2017a).

The iteration histories of objective function values and the time costs of each stage of the optimization loop are all illustrated in Fig. 15. To be specific, the time cost per iteration in LP1 (Generating TDFs and their derivatives), LP4

(Sensitivity analysis) and LP5 (Updating design variables) are denoted as $t_{TDF}$, $t_{sens}$ and $t_{updt}$, respectively. The time cost of FEA in LP3 are decomposed into two parts, i.e., the time for identification of load transmission path ($t_{srch}$) and the time for FEA with reduced DOFs ($t_{FEr}$). The average values of $t_{TDF}, t_{srch}, t_{FEr}, t_{sens}, t_{updt}$ are 2.728s, 1.038s, 3.189s, 0.153s, 0.006s, respectively. Obviously, the time costs of sensitivity analysis and updating design variables are much cheaper than the others, since fewer design variables are involved in the MMC method. For the values of $t_{TDF}$ and $t_{srch}$, they have a drastic decrease because of the elimination of tiny components. In addition, although the FEA with reduced DOFs has the largest average time cost, it is comparable with $t_{TDF}$ and $t_{srch}$. This implies, on the one hand, when the number of components is relatively large, more efficient TDF-generating method is desired (Liu et al. 2018); on the other hand, the redundant DOFs removal technique is very effective to reduce the time cost for FEA for solving 3D topology optimization problems.

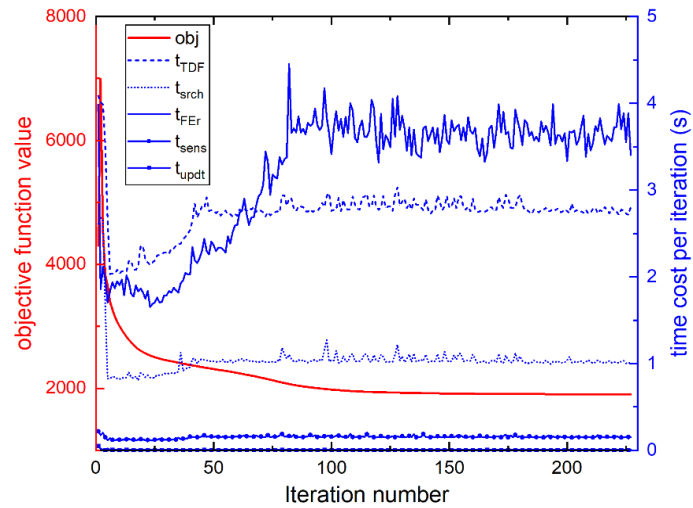

**Fig. 15** The iteration histories of the objective function values and the average time costs for each stage of the optimization loop (torsion beam example).

To test the performance of the proposed MMC code for relatively large-scale problems, this problem is solved again by refining the mesh as $192 \times 64 \times 64$ with other settings unchanged. The optimization loop converges at iteration 132 and a hollow box like structure is obtained (Fig. 16). This is also consistent with the fact that hollow boxes have excellent shear capacity.

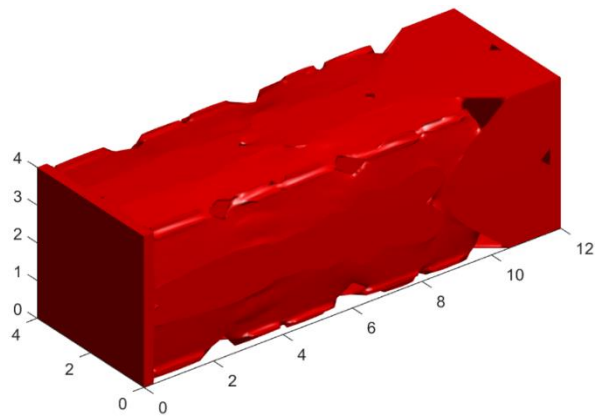

**Fig. 16** The optimized hollow box like structure (obtained with a $192 \times 64 \times 64$ mesh).

Furthermore, the optimized designs (i.e., Fig. 14(f) and Fig. 16) are analyzed for 10 times with and without using the DOFs removal technique respectively. As shown in Table 4, by refining the mesh for 1 time along each direction, the solution time of FEA with full DOFs is increased by more than 20 times. This clearly shows the issue of curse of dimensionality associated with 3D topology optimization problems. However, when the DOFs removal technique is used, both for the cases of $96 \times 32 \times 32$ mesh and $192 \times 64 \times 64$ mesh, the computational time is accelerated by more than 30 times as compared to the time cost of FEA with full DOFs. Meanwhile, the relative differences between objective function values obtained with and without DOFs removal are less than 0.3%. All these facts clearly demonstrate the efficiency of the proposed code for solving relatively large-scale 3D topology optimization problems.

**Table 4** Comparison of analysis results of the torsion beam example with/without using the redundant DOFs removal technique. The percentage values shown in the parentheses are the ratios of the involved DOFs in FEA, time cost for FEA and the relative difference of the objective function values to the corresponding values (DOFs 1, Time 1 and Obj. 1) when FEA with full DOFs is performed.

| Discretization | Elements | Nodes | DOFs 1 | DOFs 2 | Time 1 | Time 2 | Obj. 1 | Obj. 2 |
|---|---|---|---|---|---|---|---|---|
| $96 \times 32 \times 32$ | 98304 | 105633 | 316899 | 132606 (41.8%) | 113.11s | 3.53s (3.1%) | 1897.28 | 1902.86 (0.29%) |
| $192 \times 64 \times 64$ | 786432 | 815425 | 2446275 | 754032 (30.8%) | 2572.27s | 68.40s (2.7%) | 2553.81 | 2560.15 (0.25%) |

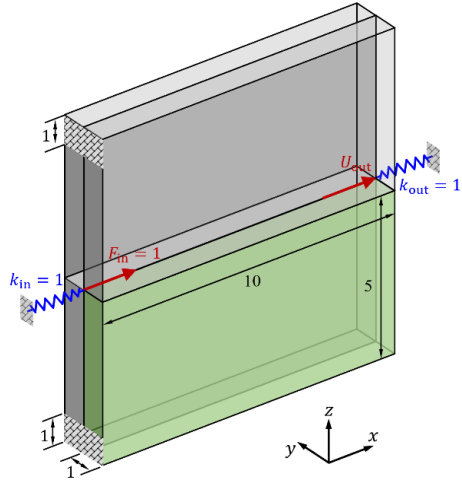

**Fig. 17** The compliant mechanism example with volume constraint of $V/|\mathrm{D}| \leq 0.2$.

### 4.4 Compliant mechanism example

Different from the stiffness design problems, an adjoint field is required for the sensitivity analysis of the compliant mechanism design. As mentioned in Remark 2, the sensitivity result in the MMC method can be calculated as a combination of a problem-dependent pre-filtered term obtained with SIMP method by taking the penalization power as 1, and a problem-independent term $\partial \rho_e / \partial d_j^i$ (different for various components), i.e.,

$$\frac{\partial U_{\mathrm{out}}}{\partial d_j^i} = -\sum_{e=1}^{N} \boldsymbol{u}_e^\top \mathbf{k}_e^0 \boldsymbol{u}_e^{\mathrm{adj}} \frac{\partial \rho_e}{\partial \phi^s} \frac{\partial \phi^s}{\partial \phi^i} \frac{\partial \phi^i}{\partial d_j^i} \tag{16}$$

where $\boldsymbol{u}_e^{\mathrm{adj}}$ is the displacement field under a pseudo unit force applied at the right-middle point along $x$-axis.

Due to the symmetry settings of the considered problem, as shown in Fig. 17, a quarter of the design domain colored in light green is discretized into $200 \times 20 \times 100$ uniform finite elements with symmetric boundary conditions. Besides, the influence of different initial designs (Fig. 18(a, d)) and different values of $\alpha$ ($10^{-3}$, $10^{-9}$) on the optimized design are investigated. Accordingly, for the initial design 1, the main code is revised as:

vInt = [2 0.5 0.25 0 atan(0.75) 0];

MMC3d(10, 1, 5, 200, 20, 100, 5/3, 1, 5/4, vInt, 0.2);

In addition, the following revisions are made in the MMC3d.m:

(1) In SEC 2, at line 14, the $\alpha$ is set as $10^{-3}$ (or $10^{-9}$), and the $\epsilon$ is set as $0.2$.

(2) In SEC 3, all the commands from lines 34 to 40 are modified as:

[iN,kN] = meshgrid(1:nelx+1,1:nelz+1);  fixNd1 = [kN*nNodfc-iN+1];

fixNd2 = [nNod-nNodfc+1:nNod];

[jN,kN] = meshgrid(1:nely+1,1:0.2*nelz+1);  fixNd3 = [(kN-1)*nNodfc+(jN-1)*(nelx+1)+1];

fixDof = unique([3*fixNd1(:)-1; 3*fixNd2(:); 3*fixNd3(:)-2; 3*fixNd3(:)-1; 3*fixNd3(:)]);

freeDof = setdiff([1:nDof],fixDof);

[jE,kE] = meshgrid(1:nely,1:0.2*nelz);  fixEle = [(kE-1)*nelx*nely+(jE-1)*nelx+1];

loadNd = nNod - nelx;

loadDof = [3*loadNd-2];

F = fsparse(loadDof',1,1/length(loadNd)/4,[nDof,1]);

loadOut = nDof - 2;

Fout = fsparse(loadOut',1,1,[nDof,1]);

loadEle = nEle - [nelx-1 0];

kin = 0.1/4;   kout = 0.1/4;

(3) In SEC 4, the initial vectors of $\beta$ and $\gamma$ in lines 52-53 are updated as:

bet = repmat([1 -1 1 -1]*vInt(5),1,N/4);

gam = repmat([1 1 -1 -1]*vInt(6),1,N/4);

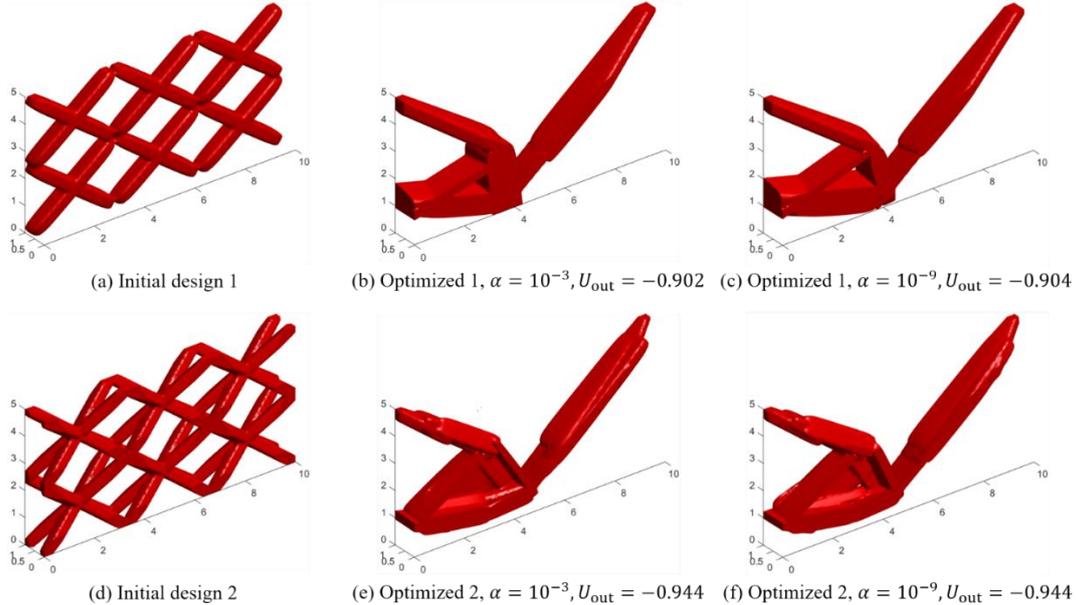

**Fig. 18** The optimized compliant mechanisms obtained with different initial designs and different alpha values.

(4) In the subsection of FEA, the stiffness of springs is added in the global stiffness matrix and the adjoint field is also calculated. Specifically, after line 92, the following command initializing $U_{\text{out}}$ should be added:

Uout = zeros(nDof,1);

Besides, line 110 is modified as:

K(loadDof,loadDof) = K(loadDof,loadDof) + kin;

K(loadOut,loadOut) = K(loadOut,loadOut) + kout;

U(freedofLft) = K(freedofLft,freedofLft)$\backslash$F(freedofLft);

Uout(freedofLft) = K(freedofLft,freedofLft)$\backslash$Fout(freedofLft);

In addition, line 115 for FEA with full DOFs is also modified as

K(loadDof,loadDof) = K(loadDof,loadDof) + kin;

K(loadOut,loadOut) = K(loadOutD,loadOut) + kout;

U(freeDof) = K(freeDof,freeDof)$\backslash$F(freeDof);

Uout(freeDof) = K(freeDof,freeDof)$\backslash$Fout(freeDof);

In line 117, the objective function value is updated as

f0val = Fout'*U/scl; OBJ(iter) = f0val*scl;

(5) In the sensitivity analysis subsection, at line 123, only the parameter energy is updated as

energy = sum((U(edofMat)*KE).*Uout(edofMat),2);

The modifications for initial design 2 are presented in Appendix 2.

It can be found in Figs. 18(a-c), in the initial design 1, there are totally 12 components distributed parallelly to the $xz$ plane. The optimized design is likely an extrusion of the 2D compliant mechanism (Sigmund, 1997; Zhang et al., 2016b). The hinge in 2D design now becomes a plane to induce the movement of the right middle point along the $-x$ direction as large as possible. Moreover, as shown in Fig. 18(d), a 3D network composed of 48 components with 432 design variables is taken as the initial design 2. The corresponding optimized designs are illustrated by Figs. 18(e,f), which are more complex in configuration and have a better performance.

This example shows that, in the MMC3d code, decreasing alpha from $10^{-3}$ to $10^{-9}$ may only have minor influence on the optimized designs, while in-appropriate initial design may lead to highly suboptimal design. This is due to that, on the one hand, structural topology optimization problems are in general highly-nonconvex in nature; while on the other hand, components can only be removed but not introduced during the optimization process. Introducing generation mechanisms of new components, for example employing topological derivatives (Takalloozadeh and Yoon, 2017; Cui et al., 2022), would be helpful to alleviate this issue.

4.5 2D cantilever beam -- test the refined 2D-MMC code

The above numerical techniques for the MMC3d code can also be implemented in a 2D-MMC code. For the $i$-th 2D component illustrated by Fig. 19, its TDF can be generated as

$$\phi^i = 1 - \left(\left(\frac{x'}{L^i}\right)^p + \left(\frac{y'}{l^i}\right)^p\right)^{1/p} \tag{17}$$

where $x'$, $y'$ and $l'$ are determined by

$$\begin{pmatrix} x' \\ y' \end{pmatrix} = \begin{bmatrix} \cos\theta^i & \sin\theta^i \\ -\sin\theta^i & \cos\theta^i \end{bmatrix} \begin{pmatrix} x - x_0^i \\ y - y_0^i \end{pmatrix} \tag{18}$$

$$l^i = \frac{t_1^i + t_2^i}{2} + \frac{t_2^i - t_1^i}{2L^i} x' \tag{19}$$

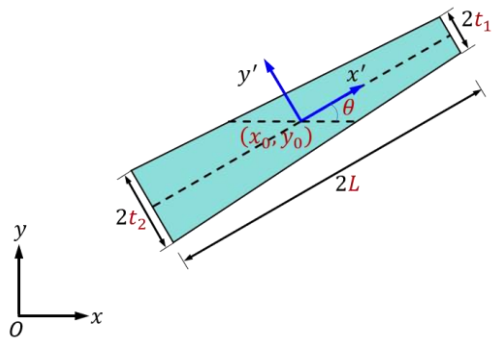

**Fig. 19** An illustration of a 2D component and the corresponding vector of the design variables $\boldsymbol{d} = (x_0, y_0, L, t_1, t_2, \theta)^\top$.

The exact expressions of $\frac{\partial \phi^i}{\partial d_j^i}$ in the 2D-MMC method are presented in Appendix 3. With the use of these expressions, the analytical sensitivities of the objective and constraint functions can be obtained using Eq. (8) and Eq. (9).

Accordingly, a 218-line Matlab code is presented in Appendix 5 for compliance minimization of the cantilever beam example illustrated in Fig. 20(a). The dimensionless Young's modulus, Poisson's ratio and the thickness are set as $E = 1, \nu = 0.3, h = 1$ respectively. The volume fraction of solid material is no greater than 40%. The design domain is discretized into $200 \times 100$ uniform quadrilateral plane stress elements and there are 16 components with equal width in the initial design. Furthermore, the MMC2d.m is activated by the following commands:

vInt = [0.4 0.05 0.05 asin(0.7)];

MMC2D(2, 1, 200, 100, 0.25, 0.25, vInt, 0.4);

At iteration 138, as illustrated in Fig. 20(b), an optimized design with structural compliance value of 73.85 is obtained ($\alpha = 10^{-9}$ and $\epsilon = 0.2$).

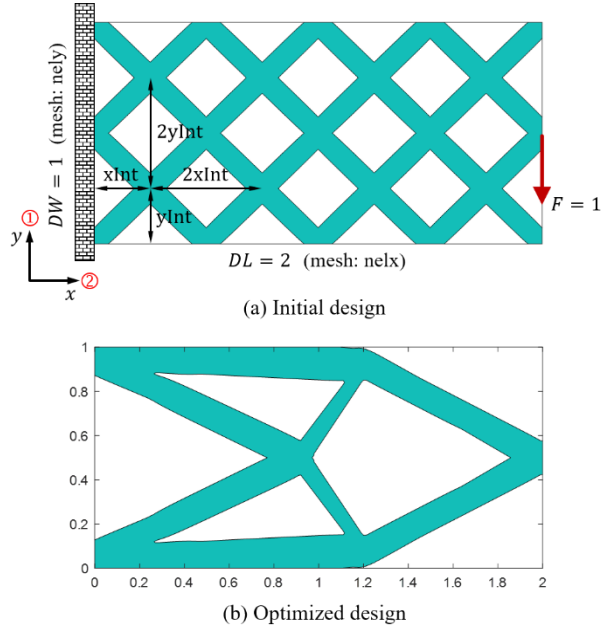

Fig. 20 The cantilever beam example solved by the MMC2d.m

## 5 Concluding remarks

In the present work, efficient and easy-to-extend Matlab codes (i.e., the 256-line MMC3d.m and the 218-line MMC2d.m) of the MMC method for structural topology optimization are developed. As compared to the existing MMC algorithms and their implementations (e.g., the MMC188 Matlab code, Zhang et al. (2016b, 2017b)), the released codes

have the following distinctive features:

(1) Analytical sensitivities which match well with the finite difference results are implemented. The proposed implementation scheme can be easily generalized to other design problems by incorporating the sensitivities obtained by the SIMP method. This helps improving the local optimality of the optimized designs and enhancing the extendibility of the proposed codes.

(2) Based on the distinctive component-based geometry description of the MMC method, a very efficient identification algorithm for load transmission path is proposed and the redundant DOFs removal technique is also implemented. This treatment effectively reduces the computational cost of FEA and makes the present codes applicable to solve relatively large-scale problems on a laptop.

(3) The strategies of identifying active components and active design variables are introduced. Tiny components have less contributions yet hard to be eliminated numerically are removed adaptively and set as inactive in the subsequent optimization step. As shown by numerical examples, this not only reduces the computational cost of generating TDFs but also accelerates the optimization process.

Although numerical examples have shown the efficiency and robustness of the MMC3d Matlab code under different mesh sizes, different initial designs and different hyper-parameters, there are still some possibilities for further improvement. For instance, introducing components with more flexible shapes to enlarge the design space Guo et al., 2016; Zhang et al., 2017b), incorporating the adaptive component introduction mechanism (Takalloozadeh and Yoon, 2017; Cui et al., 2022), calculating the TDFs and their derivatives locally to reduce the related computational cost (Liu et al., 2018), implementing the higher-order finite elements to improve the accuracy of FEA and developing the parallel version of the MMC codes for large scale 3D problems (Bathe, 2006; Aage et al., 2015).

**Acknowledgements** The authors are grateful to Prof. Krister Svanberg for providing the MMA optimizer code

**Founding** This work is supported by the National Natural Science Foundation (11821202, 11732004, 12002073, 12002077, 11922204, 11872141), the National Key Research and Development Plan (2020YFB1709401), the Fundamental Research Funds for the Central Universities (DUT20RC(3)020, DUT21-RC(3)076), Dalian Talent Innovation Program (2020RQ099), Liaoning Revitalization Talents Program (XLYC1907119), Doctoral Scientific Research Foundation of Liaoning Province (2021-BS-063) and 111 Project (B14013).

**Conflict of interest** The authors declare that they have no conflict of interest.

**Replication of results** Matlab codes are listed in the Appendix. The stenglib package, containing the fsparse function, is available at https://github.com/stefanengblom/stenglib. The MMA subroutines can be downloaded at http://www.smoptit.se/.

**References**

Aage N, Andreassen E, Lazarov BS (2015) Topology optimization using petsc: An easy-to-use, fully parallel, open source topology optimization framework. Structural and Multidisciplinary Optimization 51(3):565–572

Andreassen E, Clausen A, Schevenels M, Lazarov BS, Sigmund O (2011) Efficient topology optimization in matlab using 88 lines of code. Structural and Multidisciplinary Optimization 43(1):1–16

Bathe KJ (2006) Finite Element Procedures. Pearson Education Inc

Challis VJ (2010) A discrete level-set topology optimization code written in matlab. Structural and Multidisciplinary Optimization 41(3):453–464

Cui T, Du Z, Liu C, Sun Z, Guo X (2022) Explicit topology optimization with moving morphable component (mmc) introduction mechanism. Acta Mechanica Solida Sinica, https://doi.org/10.1007/s10338-021-00308-x

Du Z, Chen H, Huang G (2020) Optimal quantum valley hall insulators by rationally engineering berry curvature and band structure. Journal of the Mechanics and Physics of Solids 135:103784

Engblom S, Lukarski D (2016) Fast matlab compatible sparse assembly on multicore computers. Parallel Computing 56:1–17

Ferrari F, Sigmund O (2020) A new generation 99 line matlab code for compliance topology optimization and its extension to 3d. Structural and Multidisciplinary Optimization 62:2211–2228

Guo X, Zhang W, Zhong W (2014) Doing topology optimization explicitly and geometrically-a new moving morphable components based framework. Journal of Applied Mechanics 81(8):081009

Guo X, Zhang W, Zhang J, Yuan J (2016) Explicit structural topology optimization based on moving morphable components (mmc) with curved skeletons. Computer Methods in Applied Mechanics and Engineering 310:711–748

Guo X, Zhou J, Zhang W, Du Z, Liu C, Liu Y (2017) Self-supporting structure design in additive manufacturing through explicit topology optimization. Computer Methods in Applied Mechanics and Engineering 323:27–63

Hoang VN, Nguyen-Xuan H (2020) Extruded geometric-component-based 3d topology optimization. Computer Methods in Applied Mechanics and Engineering 371:113293

Huang X, Xie YM (2010) A further review of eso type methods for topology optimization. Structural and Multidisciplinary Optimization 41(5):671–683

Kreisselmeier G, Steinhauser R (1980) Systematic control design by optimizing a vector performance index. In: Computer aided design of control systems, Proceedings of the IFAC Symposium, pp 113–117

Liu C, Zhu Y, Sun Z, Li D, Du Z, Zhang W, Guo X (2018) An efficient moving morphable component (mmc)-based approach for multi-resolution topology optimization. Structural and Multidisciplinary Optimization 58(6):2455–2479

Liu K, Tovar A (2014) An efficient 3d topology optimization code written in matlab. Structural and Multidisciplinary Optimization 50(6):1175–1196

Luo J, Du Z, Guo Y, Liu C, Zhang W, Guo X (2021) Multi-class, multi-functional design of photonic topo- logical insulators by rational symmetry-indicators engineering. Nanophotonics 10(18): 4523-4531

Niu B, Wadbro E (2019) On equal-width length-scale control in topology optimization. Structural and Multidisciplinary

Optimization 59(4):1321–1334

Picelli R, Townsend S, Brampton C, Norato J, Kim HA (2018) Stress-based shape and topology optimization with the level set method. Computer methods in applied mechanics and engineering 329:1–23

Raponi E, Bujny M, Olhofer M, Aulig N, Boria S, Duddeck F (2019) Kriging-assisted topology optimization of crash structures. Computer Methods in Applied Mechanics and Engineering 348:730–752

Sigmund O (1997) On the design of compliant mechanisms using topology optimization. Journal of Structural Mechanics 25(4):493–524

Sigmund O (2001) A 99 line topology optimization code written in matlab. Structural and Multidisciplinary Optimization 21(2):120–127

Smith H, Norato JA (2020) A matlab code for topology optimization using the geometry projection method. Structural and Multidisciplinary Optimization 62:1579–1594

Svanberg K (1987) The method of moving asymptotes - a new method for structural optimization. International Journal for Numerical Methods in Engineering 24(2):359–373

Svanberg K (2007) Mma and gcmma, versions september 2007. Optimization and Systems Theory 104

Takalloozadeh M, Yoon GH (2017) Implementation of topological derivative in the moving morphable components approach. Finite Elements in Analysis and Design 134:16–26

Wang MY, Wang X, Guo D (2003) A level set method for structural topology optimization. Computer Methods in Applied Mechanics and Engineering 192(1-2):227–246

Wei P, Li Z, Li X, Wang MY (2018) An 88-line matlab code for the parameterized level set method based topology optimization using radial basis functions. Structural and Multidisciplinary Optimization 58(2):831–849

Xue R, Li R, Du Z, Zhang W, Zhu Y, Sun Z, Guo X (2017) Kirigami pattern design of mechanically driven formation of complex 3d structures through topology optimization. Extreme Mechanics Letters 15:139–144

Xue R, Liu C, Zhang W, Zhu Y, Tang S, Du Z, Guo X (2019) Explicit structural topology optimization under finite deformation via moving morphable void (mmv) approach. Computer Methods in Applied Mechanics and Engineering 344:798–818

Yoon GH, Kim YY (2003) The role of s-shape mapping functions in the simp approach for topology optimization. KSME International Journal 17(10):1496–1506

Zhang W, Li D, Zhang J, Guo X (2016a) Minimum length scale control in structural topology optimization based on the moving morphable components (mmc) approach. Computer Methods in Applied Mechanics and Engineering 311:327–355

Zhang W, Yuan J, Zhang J, Guo X (2016b) A new topology optimization approach based on moving morphable components (mmc) and the ersatz material model. Structural and Multidisciplinary Optimization 53(6):1243–1260

Zhang W, Chen J, Zhu X, Zhou J, Xue D, Lei X, Guo X (2017a) Explicit three dimensional topology optimization via moving morphable void (mmv) approach. Computer Methods in Applied Mechanics and Engineering 322:590–614

Zhang W, Li D, Yuan J, Song J, Guo X (2017b) A new three-dimensional topology optimization method based on moving morphable components (mmcs). Computational Mechanics 59(4):647–665

**Appendix 1 Exact expressions of $\frac{\partial \phi^i}{\partial d_j^i}$ in the MMC3d code**

For conciseness, the following notations are introduced

$$w \triangleq \left(\frac{x'}{L_1^i}\right)^p + \left(\frac{y'}{L_2^i}\right)^p + \left(\frac{z'}{L_3^i}\right)^p \tag{20}$$

$$\tilde{x} = x'/L_1^i, \tilde{y} = y'/L_2^i, \tilde{z} = z'/L_3^i \tag{21}$$

$$R \triangleq \begin{bmatrix} \cos\beta^i\cos\gamma^i & \cos\beta^i\sin\gamma^i & -\sin\beta^i \\ \sin\alpha^i\sin\beta^i\cos\gamma^i - \cos\alpha^i\sin\gamma^i & \sin\alpha^i\sin\beta^i\sin\gamma^i + \cos\alpha^i\cos\gamma^i & \sin\alpha^i\cos\beta^i \\ \cos\alpha^i\sin\beta^i\cos\gamma^i + \sin\alpha^i\sin\gamma^i & \cos\alpha^i\sin\beta^i\sin\gamma^i - \sin\alpha^i\cos\gamma^i & \cos\alpha^i\cos\beta^i \end{bmatrix} \tag{22}$$

$$R^\alpha \triangleq \begin{bmatrix} 0 & 0 & 0 \\ \cos\alpha^i\sin\beta^i\cos\gamma^i + \sin\alpha^i\sin\gamma^i & \cos\alpha^i\sin\beta^i\sin\gamma^i - \sin\alpha^i\cos\gamma^i & \cos\alpha^i\cos\beta^i \\ -\sin\alpha^i\sin\beta^i\cos\gamma^i + \cos\alpha^i\sin\gamma^i & -\sin\alpha^i\sin\beta^i\sin\gamma^i - \cos\alpha^i\cos\gamma^i & -\sin\alpha^i\cos\beta^i \end{bmatrix} \tag{23}$$

$$R^\beta \triangleq \begin{bmatrix} -\sin\beta^i\cos\gamma^i & -\sin\beta^i\sin\gamma^i & -\cos\beta^i \\ \sin\alpha^i\cos\beta^i\cos\gamma^i & \sin\alpha^i\cos\beta^i\sin\gamma^i & -\sin\alpha^i\sin\beta^i \\ \cos\alpha^i\cos\beta^i\cos\gamma^i & \cos\alpha^i\cos\beta^i\sin\gamma^i & -\cos\alpha^i\sin\beta^i \end{bmatrix} \tag{24}$$

$$R^\gamma \triangleq \begin{bmatrix} -\cos\beta^i\sin\gamma^i & \cos\beta^i\cos\gamma^i & 0 \\ -\sin\alpha^i\sin\beta^i\sin\gamma^i - \cos\alpha^i\cos\gamma^i & \sin\alpha^i\sin\beta^i\cos\gamma^i - \cos\alpha^i\sin\gamma^i & 0 \\ -\cos\alpha^i\sin\beta^i\sin\gamma^i + \sin\alpha^i\cos\gamma^i & \cos\alpha^i\sin\beta^i\cos\gamma^i + \sin\alpha^i\sin\gamma^i & 0 \end{bmatrix} \tag{25}$$

Then the exact formulations of $\partial \phi^i / \partial d_j^i$ can be expressed as

$$\frac{\partial \phi^i}{\partial x_0^i} = w^{\frac{1-p}{p}} \left( \frac{\tilde{x}^{p-1}}{L_1^i} R_{11} + \frac{\tilde{y}^{p-1}}{L_2^i} R_{21} + \frac{\tilde{z}^{p-1}}{L_3^i} R_{31} \right) \tag{26}$$

$$\frac{\partial \phi^i}{\partial y_0^i} = w^{\frac{1-p}{p}} \left( \frac{\tilde{x}^{p-1}}{L_1^i} R_{12} + \frac{\tilde{y}^{p-1}}{L_2^i} R_{22} + \frac{\tilde{z}^{p-1}}{L_3^i} R_{32} \right) \tag{27}$$

$$\frac{\partial \phi^i}{\partial z_0^i} = w^{\frac{1-p}{p}} \left( \frac{\tilde{x}^{p-1}}{L_1^i} R_{13} + \frac{\tilde{y}^{p-1}}{L_2^i} R_{23} + \frac{\tilde{z}^{p-1}}{L_3^i} R_{33} \right) \tag{28}$$

$$\frac{\partial \phi^i}{\partial L_1^i} = \frac{1}{L_1^i} w^{\frac{1-p}{p}} \tilde{x}^p \tag{29}$$

$$\frac{\partial \phi^i}{\partial L_2^i} = \frac{1}{L_2^i} w^{\frac{1-p}{p}} \tilde{y}^p \tag{30}$$

$$\frac{\partial \phi^i}{\partial L_3^i} = \frac{1}{L_3^i} w^{\frac{1-p}{p}} \tilde{z}^p \tag{31}$$

$$\frac{\partial \phi^i}{\partial \alpha} = -w^{\frac{1-p}{p}} \left( \frac{\tilde{x}^{p-1}}{L_1^i} \frac{\partial x'}{\partial \alpha} + \frac{\tilde{y}^{p-1}}{L_2^i} \frac{\partial y'}{\partial \alpha} + \frac{\tilde{z}^{p-1}}{L_3^i} \frac{\partial z'}{\partial \alpha} \right) \tag{32}$$

$$\frac{\partial \phi^i}{\partial \beta} = -w^{\frac{1-p}{p}} \left( \frac{\tilde{x}^{p-1}}{L_1^i} \frac{\partial x'}{\partial \beta} + \frac{\tilde{y}^{p-1}}{L_2^i} \frac{\partial y'}{\partial \alpha} + \frac{\tilde{z}^{p-1}}{L_3^i} \frac{\partial z'}{\partial \beta} \right) \tag{33}$$

$$\frac{\partial \phi^i}{\partial \gamma} = -w^{\frac{1-p}{p}} \left( \frac{\tilde{x}^{p-1}}{L_1^i} \frac{\partial x'}{\partial \gamma} + \frac{\tilde{y}^{p-1}}{L_2^i} \frac{\partial y'}{\partial \gamma} + \frac{\tilde{z}^{p-1}}{L_3^i} \frac{\partial z'}{\partial \gamma} \right) \tag{34}$$

with the following identities

$$\frac{\partial x'}{\partial \zeta} = R^{\zeta}_{11}(x - x_0) + R^{\zeta}_{12}(y - y_0) + R^{\zeta}_{13}(z - z_0) \tag{35}$$

$$\frac{\partial y'}{\partial \zeta} = R^{\zeta}_{21}(x - x_0) + R^{\zeta}_{22}(y - y_0) + R^{\zeta}_{23}(z - z_0) \tag{36}$$

$$\frac{\partial z'}{\partial \zeta} = R^{\zeta}_{31}(x - x_0) + R^{\zeta}_{32}(y - y_0) + R^{\zeta}_{33}(z - z_0) \tag{37}$$

where $\zeta = \alpha, \beta, \gamma$ respectively.

**Appendix 2 Modifications of the MMC3d code for the initial design 2 of compliant mechanism design example**

The main code is updated as

vInt = [1.2 0.2 0.2 atan(2.5) atan(0.75) atan(0.3)];

MMC3d(10, 1, 5, 200, 20, 100, 5/6, 1/4, 5/8, vInt, 0.2);

In the MMC3d.m, the modifications in SEC 2, SEC 3, the finite element analysis and sensitivity analysis subsections are the same as the modifications for initial design 1. Nevertheless, lines 44-53 in SEC 4 for the generation of initial design variables are modified as:

x0 = xInt:2*xInt:DL;  y0 = yInt:2*yInt:DW;  z0 = zInt:2*zInt:DH;

xn = length(x0);    yn = length(y0);    zn = length(z0);

x0 = kron(x0,ones(1,1*yn*zn));

y0 = repmat(repmat(y0,1,1*zn),1,1*xn);

z0 = repmat(reshape(repmat(z0,2,1),1,8),1,xn);

N = length(x0);

l1 = repmat(vInt(1),1,N);

l2 = repmat(vInt(2),1,N);

l3 = repmat(vInt(3),1,N);

alp = repmat(vInt(4),1,N);

bet = repmat([repmat([-1 -1 1 1],1,2) repmat([1 1 -1 -1],1,2)]*vInt(5),1,3);

gam = repmat([repmat([1 -1],1,4) repmat([-1 1],1,4)]*vInt(6),1,3);

**Appendix 3 Exact expressions of $\frac{\partial \phi^i}{\partial d^i_j}$ in the MMC2d code**

For conciseness, the following notations are introduced

$$w \triangleq \left(\frac{x'}{L^i_1}\right)^p + \left(\frac{y'}{L^i_2}\right)^p \tag{38}$$

$$\tilde{x} = x'/L^i, \tilde{y} = y'/l^i \tag{39}$$

And the derivatives of intermediate variables $x', y', l^i$ with respect to design variables are first detailed as

$$\begin{cases} \dfrac{\partial x'}{\partial x_0^i} = -\cos\theta^i \\ \dfrac{\partial x'}{\partial y_0^i} = -\sin\theta^i \\ \dfrac{\partial x'}{\partial \theta^i} = y' \end{cases} \quad \begin{cases} \dfrac{\partial y'}{\partial x_0^i} = \sin\theta^i \\ \dfrac{\partial y'}{\partial y_0^i} = -\cos\theta^i \\ \dfrac{\partial y'}{\partial \theta^i} = -x' \end{cases} \tag{40}$$

$$\begin{cases} \dfrac{\partial l^i}{\partial x_0^i} = \dfrac{\partial l^i}{\partial x'}\dfrac{\partial x'}{\partial x_0^i} = -\dfrac{t_2^i - t_1^i}{2L^i}\cos\theta^i \\ \dfrac{\partial l^i}{\partial y_0^i} = -\dfrac{t_2^i - t_1^i}{2L^i}\sin\theta^i \\ \dfrac{\partial l^i}{\partial L^i} = -\dfrac{t_2^i - t_1^i}{2(L^i)^2}x' \\ \dfrac{\partial l^i}{\partial t_1^i} = \dfrac{1}{2} - \dfrac{x'}{2L^i} \\ \dfrac{\partial l^2}{\partial t_2^i} = \dfrac{1}{2} + \dfrac{x'}{2L^i} \\ \dfrac{\partial l^i}{\partial \theta^i} = \dfrac{t_2^i - t_1^i}{2L^i}y' \end{cases} \tag{41}$$

Then the exact expressions of $\dfrac{\partial \phi^i}{\partial d_j^i}$ in the MMC2d code can be formulated as

$$\frac{\partial \phi^i}{\partial x_0^i} = w^{\frac{1-p}{p}}\left(\frac{\tilde{x}^{p-1}}{L^i}\cos\theta^i - \frac{\tilde{y}^{p-1}}{l^i}\sin\theta^i + \frac{\tilde{y}^p}{l^i}\frac{\partial l^i}{\partial x_0^i}\right) \tag{42}$$

$$\frac{\partial \phi^i}{\partial y_0^i} = w^{\frac{1-p}{p}}\left(\frac{\tilde{x}^{p-1}}{L^i}\sin\theta^i + \frac{\tilde{y}^{p-1}}{l^i}\cos\theta^i + \frac{\tilde{y}^p}{l^i}\frac{\partial l^i}{\partial y_0^i}\right) \tag{43}$$

$$\frac{\partial \phi^i}{\partial L^i} = w^{\frac{1-p}{p}}\left(\frac{\tilde{x}^p}{L^i} + \frac{\tilde{y}^p}{l^i}\frac{\partial l^i}{\partial L^i}\right) \tag{44}$$

$$\frac{\partial \phi^i}{\partial t_1^i} = w^{\frac{1-p}{p}}\frac{\tilde{y}^p}{l^i}\frac{\partial l^i}{\partial t_1^i} \tag{45}$$

$$\frac{\partial \phi^i}{\partial t_2^i} = w^{\frac{1-p}{p}}\frac{\tilde{y}^p}{l^i}\frac{\partial l^i}{\partial t_2^i} \tag{46}$$

$$\frac{\partial \phi^i}{\partial \theta^i} = w^{\frac{1-p}{p}}\left(-\frac{\tilde{x}^{p-1}}{L^i}y' + \frac{\tilde{y}^{p-1}}{l^i}x' + \frac{\tilde{y}^p}{l^i}\frac{\partial l^i}{\partial \theta^i}\right) \tag{47}$$